\documentclass[a4paper,12pt,oneside]{article}
\usepackage{amsfonts}
\usepackage{amsthm}
\usepackage{amssymb}
\usepackage{amsmath}
\usepackage{times}
\theoremstyle{definition}
\newtheorem{df}{\indent Definition}
\newtheorem{prz}{\indent Example}

\theoremstyle{plain}
\newtheorem{lem}[df]{\indent Lemma}
\newtheorem{tw}{\indent Theorem}
\newtheorem{stw}{\indent Proposition}
\newtheorem{uw}[df]{\indent Remark}
\newtheorem{wn}{\indent Conclusion}

\begin{document}

\author{Diana Dziewa-Dawidczyk\footnote{Department of Applied
Mathematics, Faculty of Applied Informatics and Mathematics,
Warsaw University of Life Sciences - SGGW, Ul. Nowoursynowska 159,
02-776 Warsaw, POLAND, e-mail: diana\_dziewa\_dawidczyk@sggw.pl}  %
\and Zbigniew Pasternak-Winiarski\footnote{Faculty of Mathematics
and Information Science,  Warsaw University of Technology, Ul.
Koszykowa 75, 00-662 Warsaw, POLAND, Institute of Mathematics,
University of Bia\l ystok Akademicka 2, 15-267 Bia\l ystok,
Poland, e-mail: Z.Pasternak-Winiarski@mini.pw.edu.pl }}
\title{Integration on differential spaces}

\maketitle

\begin{abstract}
\noindent The generalization of the n-dimensional cube, an
n-dimensional chain, the exterior derivative and the integral of a
differential $n$-form on it are introduced and investigated. The
analogue of Stokes theorem for the differential space is given.


\end{abstract}


\section{Introduction}
This article is the fifth of the series of papers concerning
integration of differential forms and densities on differential
spaces. We describe our generalization of the theory of
integration of smooth skew-symmetric forms on cubes and chains.

In Section 2 we present definition and properties of so called
\emph{point differential forms} on differential spaces. We prove
the theorem about the local representation of such forms (Theorem
\ref{41}). Section 3 of the paper contains basic definitions and
the description of preliminary facts concerning generalized cubes
and chains.  We define the notion of \emph{a generalized
n-dimensional smooth cube} on the differential space $(M,
{\mathcal C})$, the notion of \emph{a smooth $n$-dimensional
chain} of generalized smooth cubes on $(M, {\mathcal C})$  and the
integral of a smooth point $n$-form on a smooth $n$-dimensional
chain. In Section 4 we prove the existence of integrals for a wide
class of smooth chains and skew-symmetric forms. To formulate this
results we introduce the classes of cubes, chains and point forms
\emph{smoothly extendable with respect to the family} ${\mathcal
G}$ of generators of a differential structure ${\mathcal C}$ on a
set $M$. At the end of the paper we give an analogue of Stokes
theorem in the case of smoothly extendable forms and chains
(Theorem \ref{43}). To do that we introduce \emph{the boundary of
generalized n-dimensional chain} (Definition \ref{67}).

 Without any other explanation we use the following symbols:
${\mathbb{N}}$ - the set of natural numbers; ${\mathbb{Q}}$ - the
set of rational numbers; $\mathbb{R}$ - the set of reals.

\section{Point forms}

Let $(M,{\mathcal C})$ be a differential space and let
$TM=\bigsqcup_{m \in M}T_{m}M$ be the space tangent to
$(M,{\mathcal C})$ (the disjoint union of point tangent spaces).
Let $\pi : TM \rightarrow M$ be the natural projection, i.e. $\pi
(v) = m$ for any $v \in T_{m}M$. We endow $TM$ with the
differential structure ${\mathcal T}{\mathcal C}$ generated by the
family of function $\{\alpha \circ \pi, \alpha \in {\mathcal C}\}
\cup \{ d \alpha, \alpha \in {\mathcal C}\}$, where $d \alpha : TM
\rightarrow \mathbb{R}$, $d \alpha(v) = v(\alpha)$ (see
\cite{cuk}, Definition 2.4 and remarks after this definition).
Then the map $\pi$ is smooth with respect to differential
structures ${\mathcal C}$ and ${\mathcal T} {\mathcal C}$.

\begin{df} A map $\omega : TM\to\mathbb{R}$ is called
\emph{a point 1-form} if for any $m \in M \quad \omega_m =
\omega_{|T_{m}M}$ is an $\mathbb{R}$-linear mapping. A point
1-form $\omega$ is said to be \emph{smooth} if $\omega \in
{\mathcal T} {\mathcal C}$ (see \cite{sas}).
\end{df}

Let us take the notation:

$T^{k}M = \{(v_1, \ldots, v_k) \in TM \times \ldots \times TM :
\pi(v_1) = \ldots = \pi(v_k)\}$,

$\tau^k {\mathcal C} = ({\mathcal T}{\mathcal C} \hat{\otimes}
\ldots\hat{\otimes}{\mathcal T}{\mathcal C})_{T^{k}M}$ for $k = 1,
2, \ldots$\\
(see \cite{cuk}, remarks after Proposition 2.2).

Let $\pi_i : T^k M \rightarrow TM$ for $i=1,\ldots,k$ be a map
given by the formula
$$\pi_i (v_1, \ldots, v_k) = v_i, \qquad (v_1,
\ldots, v_k) \in T^k M.$$

\begin{df} A map $\omega : T^k M\to\mathbb{R}$ is called
\emph{a point k-form} if for any $m \in M \quad \omega_m =
\omega_{|T_{m}M}$ is a $k$-linear mapping. A point $k$-form
$\omega$ is said to be \emph{smooth} if $\omega \in \tau^k
{\mathcal C}$ (see \cite{sas}).
\end{df}

\begin{lem} Let $W$ be an open set in $\mathbb{R}^m$, $U$ -- an
open neighborhood of zero in $\mathbb{R}^n$ and $\sigma: W\times
U\to\mathbb{R}$ -- a $C^\infty$-function. If for a point $(w,u) =
(w_1, \ldots, w_m,u_1, \ldots, u_n) \in W\times U$ there exist
$a<0$ and $b>1$ such that for any $t\in(a;b)$ the point $(w,tu)\in
W\times U$  and $\sigma(w,tu) = t \sigma(w,u)$ then
$$\sigma(w,u) = \sum_{i=1}^n
\frac{\partial \sigma}{\partial u_i}(w,0) \cdot u_i.$$
\label{31}\end{lem}

\emph{Proof}. Since the map $(a;b)\ni
t\mapsto\sigma(w,tu)=t\sigma(w,u)\in\mathbb{R}$ is well defined
and smooth on $(a,b)$ we have
$$\sigma(w,u)=\frac{d}{dt}(t\sigma(w,u)_{|t=0}=\frac{d}{dt}(\sigma(w,tu)_{|t=0}
=\sum_{i=1}^n\frac{\partial \sigma}{\partial u_i}(w,0) \cdot\
\frac{d(tu_i)}{dt}$$
$$=\sum_{i=1}^n \frac{\partial\sigma}{\partial u_i}(w,0) \cdot u_i.$$
\hfill$\Box$

Using induction we can extend this result to the following one.
\begin{lem} Let $W$ be an open set in $\mathbb{R}^m$ and let for
any $j=1,2\ldots,k$ a set $U^j$ be an open neighborhood of zero in
$\mathbb{R}^{n_j}\ (n_1,\ldots,n_k\in\mathbb{N})$. Suppose that
$\sigma:W\times U^1\times\ldots\times U^k\to\mathbb{R}$ is a
$C^\infty$-function, $(w,u^1,\ldots,u^k)\in W\times
U^1\times\ldots\times U^k$ and there exist $a<0$ and $b>1$ such
that for any $t_1,\ldots,t_k\in (a;b)$ we have
$(w,t_1u^1,\ldots,t_ku^k)\in W\times U^1\times\ldots\times U^k$
and
$$\sigma(w,t_1u^1,\ldots,t_ku^k)=t_1t_2\cdots t_k\sigma
(w,u^1,\ldots,u^k).$$ Then
$$\sigma(w,u^1,\ldots,u^k)=\sum\limits_{j_1=1}^{n_1}
\cdots\sum\limits_{j_k=1}^{n_k}\frac{\partial^k\sigma} {\partial
u^1_{j_1}\ldots\partial u^k_{j_k}}(w,0,\ldots,0)u^1_{j_1}\cdots
u^k_{j_k}.$$   \hfill$\Box$ \label{32}
\end{lem}

\begin{tw} If $\omega:T^kM\to\mathbb{R}$ is a smooth point
$k$-form on a differential space $(M,{\mathcal C})$ then for any
$m\in M$ there exist a neighborhood $W_m$ of $m$, a number
$n\in\mathbb{N}$, functions $\alpha_1,\ldots,\alpha_n\in {\mathcal
C}$ and for any $i_1,\ldots,i_k\in\{1,2,\ldots,n\}$ a function
$\omega_{i_1\ldots i_k}\in{\mathcal C}$ such that
\begin{equation}
\omega_{|T^kW_m}=\sum\limits_{i_1\ldots
i_k=1}^n[(\omega_{i_1\ldots i_k}
\circ\pi\circ\pi_1)d\alpha_{i_1}\otimes\ldots\otimes
d\alpha_{i_k}]_{|T^kW_m},   \label{1}
\end{equation}
where $T^kW_m=\bigsqcup\limits_{p\in W_m}T_m^kW_m$, $\pi : TM
\rightarrow M$ is the natural projection and $\pi_i : T^k M
\rightarrow TM$ is described after Definition 1.
\label{41}\end{tw}

\begin{uw} Taking into account that functions $\omega_{i_1\ldots
i_k}$ and $\alpha_j$ depend only on $p\in W_m$ we will abbreviate
formula (\ref{1}) writing
\begin{equation}
\omega=\sum\limits_{i_1\ldots i_k=1}^n\omega_{i_1\ldots i_k}
d\alpha_{i_1}\otimes\ldots\otimes d\alpha_{i_k}\ \ \  on\ \ \ W_m.
\label{2}
\end{equation}
\end{uw}

\emph{Proof of the theorem}. It follows from the definition of the
differential structure $\tau^k {\mathcal C}$ (see \cite{cuk}) that
there exist: a neighborhood $\tilde{V}$ of the point
$(0_m,\ldots,0_m)\in T_m^kM$ in $T^kM$ (where $0_m$ is the neutral
element  of the vector space $T_m$), functions
$\beta^0_1,\ldots,\beta^0_{n_0},\beta^1_1,\ldots,\beta^1_{n_1},\ldots,
\beta^k_1,\ldots,\beta^1_{n_k}\in{\mathcal C}$ and a function
$\sigma\in C^\infty(\mathbb{R}^n)$, where $n=n_0+n_1+\ldots+n_k$,
such that for any $(v_1,\ldots,v_k)\in\tilde{V}$
$$\omega(v_1,\ldots,v_k)=$$
$$\sigma(\beta^0_1(\pi(v_1)),\ldots,
\beta^0_{n_0}(\pi(v_1)),d\beta^1_1(v_1),\ldots,d\beta^1_{n_1}(v_1),
\ldots,d\beta^k_1(v_k),\ldots,d\beta^k_{n_k}(v_k)).$$ Reducing
$\tilde{V}$, if necessary, we can assume that
$\tilde{V}=\beta^{-1}(U^0\times U^1\times\ldots\times U^k)$, where
$\beta:T^kM\to\mathbb{R}^n$ is given by the formula
$$\beta(v_1,\ldots,v_k)=$$
$$(\beta^0_1(\pi(v_1)),\ldots,
\beta^0_{n_0}(\pi(v_1)),d\beta^1_1(v_1),\ldots,d\beta^1_{n_1}(v_1),
\ldots,d\beta^k_1(v_k),\ldots,d\beta^k_{n_k}(v_k)),$$
$$(v_1,\ldots,v_k)\in T^kM$$
and $U^j\subset\mathbb{R}^{n_j}$ is an open $n_j$-dimensional
interval (cube) for $j=0,1,\ldots,k$ (see \cite{DDD}, Stwierdzenie
2.7 or \cite{sik}). Note that if $(t_1,\ldots,t_k)\in\mathbb{R}^k$
and $(v_1,\ldots,v_k)\in\tilde{V}$ fulfill the condition
$$(t_1d\beta^1_1(v_1),\ldots,t_1d\beta^1_{n_1},\ldots,t_kd\beta^k_1,
\ldots,d\beta^k_{n_k})\in U^1\times\ldots\times U^k$$
then
$$(t_1v^1,t_2v^2,\ldots,t_kv^k)\in\tilde{V}.$$

Taking $W:=(\pi\circ\pi_1)(\tilde{V}),\ \
w:=\beta^0(\pi(v_1)):=(\beta^0_1(\pi(v_1)),\ldots,
\beta^0_{n_0}(\pi(v_1))$ and $u^j:=
(d\beta^j_1(v_j),\ldots,d\beta^j_{n_j}(v_j))\in\mathbb{R}^{n_j}$
for $j=1,\ldots,k$ we obtain that the sets $W,U^1,\ldots,U^k$, the
point $(w,u^1,\ldots,u^k)$ and the map $\sigma$ fulfill conditions
of Lemma \ref{32}. By the thesis of this lemma
$$\omega(v_1,\ldots,v_k)=
\sum\limits_{j_1=1}^{n_1}\cdots\sum\limits_{j_k=1}^{n_k}
\frac{\partial^k\sigma}{\partial u_{j_1}^1\ldots\partial
u_{j_k}^k}\Big(w,0,\ldots,0\Big)d\beta^1_{j_1}(v_1) \cdots
d\beta^k_{j_k}(v_k)$$
\begin{equation}
=\sum\limits_{j_1=1}^{n_1}\cdots\sum\limits_{j_k=1}^{n_k}
\frac{\partial^k\sigma}{\partial u_{j_1}^1\ldots\partial
u_{j_k}^k}\Big(\beta^0(\pi(v_1)),0,\ldots,0\Big)d\beta^1_{j_1}
\otimes\cdots\otimes d\beta^k_{j_k}(v_1,\ldots,v_k). \label{3}
\end{equation}
Note that the function $T^kM\ni(v_1,\ldots,v_k)\mapsto
\frac{\partial^k\sigma}{\partial u_{j_1}^1\ldots\partial
u_{j_k}^k}\Big(\beta^0(\pi(v_1)),0,\ldots,0\Big)\in\mathbb{R}$
belongs to ${\mathcal C}$. Putting $\alpha_i:=\beta^0_i$ for
$i=1,2,\ldots,n_0$, $\alpha_{n_0+n_1+\ldots+n_l+i}:=\beta_i^{l+1}$
for $l=0,1,\ldots,k-1$ and $i=1,\ldots,n_{l+1}$,
$n:=n_0+n_1+\ldots+n_k$, $\omega_{i_1\ldots i_n}(q):=
\frac{\partial^k\sigma}{\partial u_{j_1}^1\ldots\partial
u_{j_k}^k}\Big(\beta^0(q),0,\ldots,0\Big)$  if
$d\alpha_{i_1}\otimes\ldots\otimes d\alpha_{i_k}= d\beta^1_{j_1}
\otimes\cdots\otimes d\beta^k_{j_k}$ ($q\in W$) and
$\omega_{i_1\ldots i_n}(q):=0$ if
$d\alpha_{i_1}\otimes\ldots\otimes d\alpha_{i_k}$ does not appear
in (\ref{3}), we obtain that (\ref{1}) is true if we change
$T^kW_m$ onto $\tilde{V}$. Hence, by $k$-linearity of $\omega$, it
follows that (\ref{1}) holds for $T^kW_m$, where $W_m=W$.
\hfill$\Box$

\smallskip
For any smooth map $F$ of a differential space $(M,{\mathcal C})$
into a differential space $(N, {\mathcal D})$ (we denote it
writing  $F : (M,{\mathcal C}) \rightarrow (N, {\mathcal D})$) and
any smooth $k$-form $\omega$ on $(N,{\mathcal D})$ one can define
the smooth $k$-form $F^* \omega$ by the formula:
$$(F^* \omega)(v_1, \ldots, v_k) = \omega(TFv_1, \ldots,TFv_k),
\qquad (v_1, \ldots, v_k) \in T^k M,$$ where the tangent map $TF:
TM \rightarrow TN$ is defined as follows
$$(TFv)(f) = v(f \circ F),\qquad v \in
TM,\ \ f \in {\mathcal D}.$$ Note that for any $m\in M$ the map
$TF_{|T_mM}$ is linear.

\begin{stw} If we denote
$F=(\alpha_0,\ldots,\alpha_n):W_m\to\mathbb{R}^n$ in the proof of
Theorem \ref{41} then the equality (\ref{2}) can be written in the
form
$$\omega_{|W_m}=(F^*\eta)_{|W_m}$$
(here we consider $\omega$ as a section of the $k$-power of
cotangent bundle on $M$), where $\eta=\sum\limits_{i_1\ldots
i_k=1} ^n\eta_{i_1\ldots i_k}du_{i_1}\otimes\cdots\otimes
du_{i_k}$ is the smooth $k$-linear form on $\mathbb{R}^n$ given as
follows
$$\eta_{i_1\ldots i_k}(u)=\frac{\partial^k\sigma}{\partial
u_{n_0+i_1}\ldots\partial
u_{n_0+n_1+\ldots+n_{k-1}+i_k}}(u_1,\ldots,u_{n_0},0,\ldots,0)$$
if $d\alpha_{i_1}\otimes\ldots\otimes d\alpha_{i_k}=
d\beta^1_{j_1} \otimes\cdots\otimes d\beta^k_{j_k}$ for some $1\le
j_1\le n_1,\ldots,1\le j_k\le n_k$ and
$$\eta_{i_1\ldots i_k}(u)=0$$
in other cases. \label{51}
\end{stw}

\emph{Proof}. It follows from the fact that for any $n_0<i<n$, any
$q\in W_m$ and any $v\in T_qM$
$$F^*du_i(v)=(TFv)(u_i)=v(u_i\circ F)=v(\alpha_i)=d\alpha_i(v).$$
\hfill$\Box$

\begin{uw} It is easy to see that if ${\mathcal G}$ is an arbitrary
family of generators of the differential structure ${\mathcal C}$
then in the formulation of Theorem \ref{41} and Proposition
\ref{51} we can put ${\mathcal G}$ instead of ${\mathcal C}$
preserving only the condition: $\omega_{i_1\ldots i_k}\in{\mathcal
C}$. For the definition and basic properties of a set of
generators of a differential structure see \cite{DDD}, \cite{cuk},
\cite{DPW1} or \cite{DPW2}. \label{21}
\end{uw}

\smallskip
In the following part of the paper we will be interested only in
skew-symmetric $k$-forms. In this case equalities (\ref{1}) and
(\ref{2}) take the form
\begin{equation}
\omega=\sum\limits_{1\le i_1<\ldots<i_k\le n}\omega_{i_1\ldots
i_k} d\alpha_{i_1}\wedge\ldots\wedge d\alpha_{i_k}\ \ \  on\ \ \
W_m. \label{4}
\end{equation}

\section{Definition of the integral on generalized cubes and
chains}
\begin{df}
The Cartesian product $[0,1]^n$ of $n$ unit intervals is called
\emph{ the standard n-dimensional interval}. A continuous mapping
$\varphi : [0,1]^n \rightarrow A$ is said to be
\emph{n-dimensional cube} in a topological space $A$. If $(M,
{\mathcal C})$ is a differential space, then the cube $\varphi :
[0,1]^n \rightarrow M$ is \emph{smooth} when $\varphi : ([0,1]^n,
C^\infty ([0,1]^n)) \rightarrow (M, {\mathcal C})$, i.e. $\varphi$
is a smooth mapping on the differential space $([0,1]^n, C^\infty
([0,1]^n))$ into the differential space $(M, {\mathcal C})$. The
integral of a smooth point $n$-form $\omega$ (on $M$) on a smooth
cube $\varphi$ is defined as the number
$$\int_\varphi \omega := \int_{[0,1]^n} \varphi^* \omega,$$ where
$\varphi^* \omega$ is the pullback of the form $\omega$ onto the
interval $[0,1]^n$.
\end{df}

\begin{prz}
Let $M=\mathbb{Q}$ be the set of rational numbers and ${\mathcal
C} := C^\infty(\mathbb{R})_\mathbb{Q}$. Then for any number $n \in
\mathbb{N}$ only constant functions are n-dimensional cubes in
$M$. If $\varphi$ is a constant cube (of any dimension) and
$\omega$ is a point form on $\mathbb{Q}$ (the same dimension),
than $\varphi^* \omega \equiv 0$. So we have $\int_\varphi \omega
= 0$. $\qquad \Box$ \label{97}\end{prz}

\smallskip
It follows from the example above that the standard integration
theory might be trivial on some differential spaces. In fact, on
the space $(\mathbb{Q},C^\infty(\mathbb{R})_\mathbb{Q})$ there are
non trivial point 1-forms, such as $\omega = d(id_\mathbb{Q})$,
but all integrals of such forms vanish. We would like to
generalize integration theory in a way, which enables us to omit
this inconvenience. We have to give the generalization of
n-dimensional cubes and chains. In order to do this we will use
the theory of completions and compactifications of a differential
space described in \cite{cuk}. At first we will extend smooth
functions.

\bigskip
Let $(M,{\mathcal C})$ be a differential space such that $M$ is
Hausdorff space and ${\mathcal C}$ is generated by the family of
functions ${\mathcal G}$, which defines uniform structure
${\mathcal U}_{\mathcal G}$ on $M$ (see \cite{cuk}). Any function
$f \in {\mathcal C}$ which is uniform with respect to ${\mathcal
U}_{\mathcal G}$ and standard uniform structure on $\mathbb{R}$
can be extended in the univocal way to the continuous function on
the completion $compl_{\mathcal G} M$ of the uniform space $(M,
{\mathcal U}_{\mathcal G})$. In particular any function $g \in
{\mathcal G}$ can be extended to the continuous function
$\widetilde{g}$ on $compl_{\mathcal G} M$. The localization of the
differential structure $compl_{\mathcal G} {\mathcal C}$ generated
by the family $\widetilde{{\mathcal G}} = \{\widetilde{g} :
compl_{\mathcal G} M \rightarrow \mathbb{R} : g \in {\mathcal
G}\}$ to the set $M$ is equal to ${\mathcal C}$. If ${\mathcal G}
= {\mathcal C}$, then any function $\alpha \in {\mathcal C}$ can
be extended to $\widetilde{\alpha} \in compl_{\mathcal C}
{\mathcal C}$ on the maximal differential completion
$compl_{\mathcal C} M$.

\bigskip
\begin{df} \emph{A generalized n-dimensional smooth cube} on the
differential space $(M, {\mathcal C})$ is any smooth mapping
$\varphi : (D,C^\infty(\mathbb{R}^n)_D) \rightarrow (M, {\mathcal
C})$ defined on the set $D$ dense in $[0;1]^n$.
\end{df}

\bigskip
If ${\mathcal G}$ is such a set of the generators of the structure
${\mathcal C}$, that $\varphi$ is uniform with respect to uniform
structure ${\mathcal U}_{\mathcal G}$ on $M$ (and standard uniform
structure on $D$), then $\varphi$ can be extended to n-dimensional
cube $\widetilde{\varphi} : [0;1]^n \rightarrow compl_{\mathcal G}
M$ (not necessarily smooth with respect to $compl_{\mathcal G}
{\mathcal C}$).

\begin{prz} Let $M$ be the set of all square roots of rational numbers
belonging to the interval $(0;1)$. A mapping
$$\varphi(x) = \sqrt{x}; \quad x \in (0;1) \cap \mathbb{Q},$$ is a
generalized smooth 1-dimensional cube in the space $(M,
C^\infty(\mathbb{R})_M)$ and it is uniform with respect to the
uniform structure ${\mathcal U}_{\mathcal G}$ defined by
${\mathcal G}=\{id_M\}$. But its extension $\widetilde{\varphi} :
[0;1] \rightarrow compl_{\mathcal G} M = [0;1]$,
$$\widetilde{\varphi} (x) = \sqrt{x}, \quad x \in
[0;1],$$ is not smooth with respect to $compl_{\mathcal G}
{\mathcal C} = C^\infty([0;1])$. $\qquad \Box$ \label{89}\end{prz}

\smallskip
If $\varphi : D \rightarrow M$ is a generalized $n$-dimensional
smooth cube on the differential space  $(M, {\mathcal C})$, and
$\omega$ is a smooth point $n$-form on $M$, then the pullback
$\varphi^* \omega$ is a smooth point $n$-form on $D$. That form is
represented by the formula:
$$[\varphi^* \omega](t) = \omega_0 (t) dt^1 \wedge \ldots \wedge dt^n,
\quad t = (t^1, \ldots, t^n) \in D,$$ where $\omega_0$ is a smooth
function on $D$. If the function $\omega_0$ could be extended to a
continuous function  $\widetilde{\omega}_0 : [0,1]^n \rightarrow
\mathbb{R}$ (the function $\widetilde{\omega}_0$ is given
univocally because $D$ is a dense set in $[0,1]^n$), then we call
the form
$$\widetilde{\varphi^* \omega} := \widetilde{\omega}_0 dt^1 \wedge
\ldots \wedge dt^n$$ \emph{the continuous extension} of $\varphi^*
\omega$ on the cube $[0,1]^n$. In that case we can define the
integral.

\begin{df}
Let $\omega$ be a smooth point $n$-form on the differential space
$(M, {\mathcal C})$. \emph{The integral} of $\omega$ on the
generalized $n$-dimensional smooth cube $\varphi : D \rightarrow
M$ is the number
$$\int_\varphi \omega := \int_{[0,1]^n} \widetilde{\varphi^* \omega},$$
where $\widetilde{\varphi^* \omega}$ is the continuous extension
of $n$-form $\varphi^* \omega$ on the cube $[0,1]^n$ (if the
extension exists).
\end{df}

\smallskip
Similarly as in the case of differential manifolds, we can treat
the integral of the smooth point form on the generalized
$n$-dimensional cube $\varphi$ as a value of the mapping given on
the set of such generalized smooth cubes, that the integral
exists. We can extend that mapping on the set of $n$-dimensional
chains on the differential space $(M, {\mathcal C})$.

\begin{df}
Let $Cube_n (M, {\mathcal C})$ be the set of all $n$-dimensional
smooth generalized cubes on a differential space $(M, {\mathcal
C})$. \emph{A smooth $n$-dimensional chain} on $(M, {\mathcal C})$
is any function $\Phi : Cube_n(M,{\mathcal C}) \rightarrow
\mathbb{R}$ such that its values are different from zero on at
most finite subset of a set $Cube_n (M,{\mathcal C})$. We denote
the set of all smooth $n$-dimensional chains on $(M,{\mathcal C})$
by $Ch_n(M,{\mathcal C})$.
\end{df}

\smallskip
The set $Ch_n(M,{\mathcal C})$ has the natural structure of a
vector space over $\mathbb{R}$, as a set of real-valued functions
on $Cube_n (M,{\mathcal C})$. "Zero chain" is the function equal
to zero at each cube. We identify any element $\varphi$ of a set
$Cube_n (M,{\mathcal C})$ with a chain $\Phi \in Ch_n (M,{\mathcal
C})$ such that $\Phi(\varphi) = 1$ and $\Phi(\psi) = 0$ for $\psi
\in Cube_n (M,{\mathcal C}) \setminus \{\varphi\}$. We denote that
chain by $1 \cdot \varphi$.

If $\Phi \in Ch_n (M,{\mathcal C})$, then we define \emph{the
support of} $\Phi$ as the set $supp \Phi = \{\varphi \in Cube_n
(M,{\mathcal C}) : \Phi (\varphi) \neq 0\}$. From the definition
of a chain we get that the set $supp \Phi$ is finite. If $supp
\Phi = \{\varphi_k, \ldots, \varphi_m\}$, where $k,m \in
\mathbb{N}$, $k \leq m$ and $\Phi (\varphi_j) = c_j$ for $j = k,
\ldots, m$, than the chain $\Phi$ can be represented by the
formula
\begin{equation}
\Phi = \sum_{j=k}^m c_j \varphi_j \ . \label{64}
\end{equation}

\smallskip
If $\Phi^{'} = \sum_{j=1}^{m} c_j \varphi_j$, $\Phi^{''} =
\sum_{j=m+1}^{n} c_j \varphi_j$, then $\Phi^{'} + \Phi^{''} =
\sum_{j=1}^{n} c_j \varphi_j$. Assuming that the input of zero
components into the sum (\ref{64}) is zero and taking $c_\varphi
:= \Phi(\varphi)$ we can write $\Phi = \sum_{\varphi \in Cube_n
(M, {\mathcal C})} c_\varphi \varphi$.

\smallskip
If $\omega$ is a smooth $n$-dimensional point form on $(M,
{\mathcal C})$, then by $Cube_n^\omega (M, {\mathcal C})$ we
denote the set of all $\varphi \in Cube_n (M, {\mathcal C})$, such
that there exists continuous extension $\widetilde{\varphi^*
\omega}$ of the form $\varphi^* \omega$ onto $[0,1]^n$. Similarly,
by $Ch_n^\omega (M, {\mathcal C})$ we denote the set of all chains
$\Phi = Ch_n (M, {\mathcal C})$, such that $supp \Phi \subset
Cube_n^\omega (M, {\mathcal C})$.

\begin{df}
\emph{The integral of a smooth point $n$-form} $\omega$ on
$(M,\mathcal{C})$ on a smooth chain $\Phi = \sum_{j=1}^{n} c_j
\varphi_j \in Ch_n^\omega (M, {\mathcal C})$ is the number
$$\int_\Phi \omega = \sum_{j=1}^{m} c_j \int_{\varphi_j} \omega.$$
\end{df}

\smallskip
The mapping $Ch_n^\omega (M, {\mathcal C}) \ni \Phi \mapsto
\int_\Phi \omega \in \mathbb{R}$ is, of course, a linear
functional on $Ch_n^\omega (M, {\mathcal C})$. We will call the
mapping \emph{'' integral of the $n$-form $\omega$''}, too.

\section{The existence of integrals and Stokes theorem}
The problem is that the domain $Ch_n^\omega (M, {\mathcal C})$ of
the integral depends on the form $\omega$ (in the case of
integrals of smooth $n$-forms on ordinary chains on a smooth
manifold we can integrate each smooth $n$-form on each chain of
class $C^1$). In the first part of this section we will give
sufficient conditions for existence of  the continuous or smooth
extension $\widetilde{\varphi^* \omega}$ of a point $n$-form
$\varphi^* \omega$ from the domain of the mapping $\varphi$ to
$[0,1]^n$, where $\omega$ is a smooth point $n$-form and $\varphi$
is a generalized smooth $n$-dimensional cube on the differential
space $(M, {\mathcal C})$ .

If $D$ is a dense subset in $[0,1]^n$, then the space $T_t D$
tangent to the differential space $(D,C^\infty(\mathbb{R}^n))$  at
the point $t=(t^1,\ldots,t^n) \in D$ is an $n$-dimensional vector
space which can be identified with
$T_t\mathbb{R}^n\cong\mathbb{R}^n$. Its base consists of vectors
$\frac{\partial}{\partial t^1}|_t ,\ldots,\frac{\partial}{\partial
t^n}|_t$. Any mapping $D \ni t \mapsto \frac{\partial}{\partial
t^j}|_t \in TD$ is a smooth vector field on $D$ (see \cite{DPW2}).
On the other hand on $TD \subset T [0,1]^n = [0,1]^n \times
\mathbb{R}^n$ there exists the standard uniform structure
inherited from the space $[0,1]^n \times \mathbb{R}^n$ (given by
coordinates of the mapping $TD \ni \sum_{j=1}^n v^j
\frac{\partial}{\partial t^j}|_t \mapsto
(t^1,\ldots,t^n,v^1,\ldots,v^n) \in [0,1]^n \times \mathbb{R}^n$)
(see \cite{DPW1}). We have:

\begin{stw}
Any vector field $\frac{\partial}{\partial t^i} : D \rightarrow
TD$ is an uniform mapping with respect to standard uniform
structures on $D$ and $TD$.
\end{stw}

\emph{Proof}. Let for given $\varepsilon > 0$
$$V_{\varepsilon} := \{(\sum_{j=1}^n v^j \frac{\partial}{\partial t^j}|_t,
\sum_{j=1}^n w^j \frac{\partial}{\partial t^j}|_s) \quad \in$$
$$\in \quad TD \times TD : |t^j - s^j|<\varepsilon \quad
\wedge \quad |v^j - w^j|<\varepsilon, \quad j=1,2,\ldots,n\}.$$ A
family $\{V_{\varepsilon}\}_{\varepsilon \in (0; + \infty)}$ is a
base of the standard uniform structure on $TD$, and a family
$\{W_{\varepsilon}\}_{\varepsilon \in (0; + \infty)}$, where
$$W_{\varepsilon} := \{(t,s) \in D \times D : |t^j - s^j| < \varepsilon,
\quad j=1,2,\ldots,n\}$$ is a base of the standard uniform
structure on $D$. It is clear that if for any $j \in \mathbb{N}$
$(t,s)=(t^1,\ldots,t^n,s^1,\ldots,s^n) \in W_{\varepsilon}$, than
$(\frac{\partial}{\partial t^j}|_t , \frac{\partial}{\partial
s^j}|_s) \in V_{\varepsilon}$. So we have the thesis (see
\cite{DPW1} and \cite{DPW2}). $\qquad \Box$

\begin{stw}
Let $\mathcal G$ be a family of generators of a differential
structure $\mathcal C$ on a set $M$. If $\varphi : D \rightarrow
M$ is such a smooth $n$-dimensional generalized cube on $M$ that
the tangent mapping $T \varphi :TD \rightarrow TM$ is uniform with
respect to the uniform structure given on $TM$ by the family
${\mathcal T} {\mathcal G}_0$, then for any function $\alpha \in
{\mathcal G}$ the point 1-form $\varphi^* d \alpha$ smooth on $TD$
can be extended to the 1-form $\widetilde{\varphi^* d \alpha}$
continuous on $[0,1]^n$. \label{84}\end{stw}

\emph{Proof}. The form $\varphi^* d \alpha$ is represented by the formula
$$\varphi^* d \alpha = \sum_{i=1}^n \frac{\partial (\alpha \circ \varphi)
(t)}{\partial t^i} d t^i, \quad t \in D,$$
where for any $i=1,2,\ldots,n$
$$\frac{\partial (\alpha \circ \varphi)}{\partial t^i} (t) =
(\varphi^* d \alpha) \frac{\partial}{\partial t^i}|_t = d \alpha
(T \varphi (\frac{\partial}{\partial t^i}|_t)).$$ Since the
mapping $T \varphi$ is uniform on $TD$, the vector field
$\frac{\partial}{\partial t^i}$ is a uniform mapping with respect
to $D$ in $TD$, and $d \alpha$ is one of functions of the family
${\mathcal T}{\mathcal G}_0$, so $\frac{\partial (\alpha \circ
\varphi)}{\partial t^i}$ is a uniform function on $D$ and it can
be extended to the continuous function $\widetilde{\frac{\partial
(\alpha \circ \varphi)}{\partial t^i}}$ on $[0,1]^n$ (see
\cite{DPW1} and \cite{DPW2}). So $\varphi^* d \alpha$ might be
extended to the continuous 1-form $\widetilde{\varphi^* d \alpha}
= \sum_{i=1}^n \widetilde{\frac{\partial(\alpha \circ
\varphi)}{\partial t^i}}(t)dt^i$. $\qquad \Box$

\begin{wn}
Let the assumptions of the previous proposition be satisfied. Then
for all $\alpha^1,\ldots,\alpha^k \in {\mathcal G}$ the point
$k$-form $\varphi^* (d \alpha^1 \wedge \ldots \wedge d \alpha^k)$
can be extended to the $k$-form $\widetilde{\varphi^*(d \alpha^1
\wedge \ldots \wedge d \alpha^k)}$, which is continuous on
$[0,1]^n$. \label{103}\end{wn}

\emph{Proof}. We have:
$$\varphi^* (d \alpha^1 \wedge \ldots \wedge d \alpha^k)
=(\varphi^* d \alpha^1) \wedge \ldots \wedge (\varphi^* d
\alpha^k).$$ So $\varphi^* (d \alpha^1 \wedge \ldots \wedge d
\alpha^k)$ could be extended to the $k$-form
$$\widetilde{\varphi^* (d \alpha^1 \wedge \ldots \wedge d \alpha^k)}
:= (\widetilde{\varphi^* d \alpha^1}) \wedge \ldots \wedge
(\widetilde{\varphi^* d \alpha^k}). \qquad \Box$$

\begin{stw}
Let the assumptions of Proposition \ref{84} be satisfied. Let $k,m
\in \mathbb{N}$, $\alpha^1_j,\ldots,\alpha^k_j \in {\mathcal G}$
for $j=1,2,\ldots,m$ and functions $\beta_1,\ldots,\beta_m : M
\rightarrow \mathbb{R}$ are uniform with respect to the uniform
structure ${\mathcal U}_{\mathcal G}$ on $M$, then the point
$k$-form $\eta :=\varphi^*(\sum_{j=1}^m \beta_j d \alpha^1_j
\wedge \ldots \wedge d \alpha^k_j)$ can be extended to the
$k$-form $\tilde{\eta}$, which is continuous on $[0,1]^n$.
\label{87}\end{stw}

\emph{Proof}. Since $\eta = \sum_{j=1}^m (\beta_j \circ \varphi)
\varphi^* (d \alpha_j^1 \wedge \ldots \wedge d \alpha_j^k)$ it is
enough to take $\tilde{\eta} = \sum_{j=1}^m (\widetilde{\beta_j
\cdot \varphi}) \cdot \widetilde{\varphi^*(d \alpha_j^1 \wedge
\ldots \wedge d \alpha_j^k)}$ (see Conclusion \ref{103}). $\Box$

\begin{stw}
Let $\eta$ be a smooth point $k$-form on a differential space
$(M,{\mathcal C})$, and ${\mathcal G}$ be a family of generators
of the differential structure ${\mathcal C}$. Then for any $p \in
M$ there exists a neighborhood $U$, a number $m \in \mathbb{N}$,
functions $\alpha_j^1, \ldots, \alpha_j^k \in {\mathcal G}$ for
$j=1,2,\ldots,m$ and functions $\beta_1,\ldots,\beta_m \in
{\mathcal C}$ uniform with respect to the uniform structure
${\mathcal U}_{\mathcal G}$ such that
$$\eta_{|U} = \sum_{j=1}^m \beta_{j_{|U}} (d \alpha_j^1)_{|U}
\wedge \ldots \wedge (d \alpha_j^k)_{|U}.$$
\end{stw}

\emph{Proof}. It follows from Theorem \ref{41} and Remark \ref{21}
that there exists a neighborhood $W_p$ of $p$ such that $\eta$ can
be written in the form (\ref{4}) on $W_p$. Let us fix any 1-1 map
$J:\{1,2,\ldots,\frac{n!}{k!(n-k)!}\}\to\{(i_1,\ldots,i_k)\in
\mathbb{N}^k:1\le i_1<i_2<\ldots<i_k\le n\}$ and let
$J(j)=(i_1(j),\ldots,i_k(j))$. Then putting
$m:=\frac{n!}{k!(n-k)!},\ \  \omega_j:=\omega_{J(j)}$ and
$\alpha_j^r:=\alpha_{i_r(j)}$ for any $j\in\{1,2,\ldots,m\}$ and
$1\le r\le k$ we obtain that
\begin{equation}
\eta_{|W_p} = \sum_{j=1}^m \omega_{J(j)|W_p} (d \alpha_j^1 \wedge
\ldots \wedge (d \alpha_j^k))_{|W_p}, \label{85}\end{equation}
where $\alpha_j^1, \ldots, \alpha_j^k\in{\mathcal G}$ for
$j=1,2,\ldots,m$ and functions $\omega_1,\ldots,\omega_m \in
{\mathcal C}$ (the number $m$, functions
$\alpha_j^1,\ldots,\alpha_j^k$ and $\omega_1,\ldots,\omega_m$
could depend on $W_p$). Recall now that $\omega_j=0$ or $\omega_j$
is the superposition of the smooth map $\beta_0\in{\mathcal
G}^{n_0}$ which takes a value in an $n_0$-dimensional interval
$U^0$ and a smooth function (partial derivative
$\frac{\partial^k\sigma}{\partial u_{j_1}^1\ldots\partial
u_{j_k}^k}$ ) defined on $\mathbb{R}^{n_0}$ (see proof of Theorem
\ref{41}). Since $\frac{\partial^k\sigma}{\partial
u_{j_1}^1\ldots\partial u_{j_k}^k}_{|U^0}$ is an uniform function
on $U^0$  and coefficients of $\beta^0$ are uniform functions with
respect to the uniform structure ${\mathcal U}_{\mathcal G}$ on
$M$ (as members of ${\mathcal G}$) we obtain that $\omega_j$ are
uniform with respect to ${\mathcal U}_{\mathcal G}$ on the set
$U:=W_p$ as the superposition of uniform mappings. The thesis
follows now from (\ref{85}). $\Box$

\begin{tw}
Let $\eta$ be a smooth point $k$-form on a differential space
$(M,{\mathcal C})$, and ${\mathcal G}$ be a family of generators
of the differential structure ${\mathcal C}$. Let us assume that:\\
(i) there exists an open  finite covering $\{U_i\}_{i\in I}$
($I=\{1,\ldots,q\},\ q\in\mathbb{N}$) of the space $M$, such that
for any $i\in I$ there exist a number $m_i \in \mathbb{N}$,
functions $\alpha_{i,j}^1,\ldots, \alpha_{i,j}^k \in {\mathcal G}$
for $j=1,\ldots,m_i$ and functions
$\beta_{i,1},\ldots,\beta_{i,m_i} \in {\mathcal C}$ uniform with
respect to the uniform structure ${\mathcal U}_{\mathcal G}$ such
that $\eta_{|U_i} = \sum_{j=1}^{m_i} \beta_{{i,j}_{|U_i}}(d
\alpha_{i,j}^1)_{|U_i} \wedge \ldots
\wedge (d \alpha_{i,j}^k)|_{U_i}$;\\
(ii) there exist a smooth partition of unity $\{\gamma_i\}_{i \in
I}$ subordinated to the covering $\{U_i\}_{i \in I}$ such that any
function $\gamma_i$, $i \in I$, is uniform with respect to
${\mathcal U}_{\mathcal G}$. Then for any smooth n-dimensional
generalized cube $\varphi : D \rightarrow M$ such that $T \varphi$
is uniform with respect to ${\mathcal T}{\mathcal G}_0$, k-form
$\varphi^* \eta$ is smooth on $D$ and could be extended to the
$k$-form $\widetilde{\varphi^* \eta}$ continuous on $[0,1]^n$.
\label{88}\end{tw}

\emph{Proof}. We have $\eta = \sum_{i \in I} \gamma_i \eta$,
and $\varphi^* \eta = \sum_{i \in I} \varphi^*(\gamma_i \eta)$.
But for any $i \in I$ $supp \gamma_i \subset U_i$ so we have
$$\gamma_i \eta = \sum_{j=1}^{m_i} \gamma_i \cdot \beta_{i,j}
d \alpha_{i,j}^1 \wedge \ldots \wedge d \alpha_{i,j}^k.$$
Therefore
$$\varphi^*(\gamma_i \eta) = \sum_{j=1}^{m_i}(\gamma_i \circ \varphi)
\varphi^* [\beta_{i,j} d \alpha_{i,j}^1 \wedge \ldots \wedge d
\alpha_{i,j}^k]$$ and the form could be extended on $[0,1]^n$ to
the $k$-form
$$\widetilde{\varphi^*(\gamma_i \eta)} := \sum_{j=1}^{m_i}
(\widetilde{\gamma_i \circ \varphi}) \varphi^*
[\widetilde{\beta_{i,j} d \alpha_{i,j}^1 \wedge \ldots \wedge d
\alpha_{i,j}^k}]$$ (see Proposition \ref{87}). Hence $\varphi^*
\eta$ could be extended to the $k$-form
$$\widetilde{\varphi^*
\eta} = \sum_{i \in I}\widetilde{\varphi^*(\gamma_i \eta)}. \qquad
\qquad \Box$$

If a point $k$-form $\eta$ and a smooth generalized
$n$-dimensional cube $\varphi$ on a differential space
$(M,{\mathcal C})$ satisfy the assumptions of Theorem \ref{88},
then there exists $\int_{\varphi} \eta$ (of course it is equal to
zero for $n \neq k$). Similarly if $\Psi$ is a smooth
$n$-dimensional chain on $(M,{\mathcal C})$ and each cube $\varphi
\in supp \Psi$ together with a point $k$-form $\eta$ satisfy the
assumptions of Theorem \ref{88} then there exists $\int_{\Psi}
\eta$.

For the further development of the theory of integration we need a
generalization of the operation of boundary. The operation
prescribes to any $n$-dimensional cube $\varphi$ an
$(n-1)$-dimensional chain $\partial \varphi$. It is connected with
Stokes theorem and justifies the concept of chain. Standardly we
can extend this operation to the set of chains. For any
$n$-dimensional chain $\Phi = \sum_{j=1}^n c_j \varphi_j$ we have
$\partial \Phi = \sum_{j=1}^n c_j \partial \varphi_j$.

To introduce the operation of boundary on generalized cubes and
chains in differential spaces we need some new concepts.\\

\begin{df}
Let ${\mathcal G}$ be a family of generators of a differential
structure ${\mathcal C}$ on a set $M$. We call a smooth
generalized $n$-dimensional cube $\varphi : D \rightarrow M$ on
$(M,{\mathcal C})$ \emph{smoothly extendable with respect to}
${\mathcal G}$ if the tangent mapping $T \varphi$ is uniform with
respect to the uniform structure given on $TM$ by the family
${\mathcal T} {\mathcal G}_0$ (see \cite{cuk}) and its continuous
extension $\widetilde{T \varphi} : [0;1]^n \times \mathbb{R}^n
\rightarrow compl_{{\mathcal T}{\mathcal G}_0}TM$ is a smooth map
with respect to differential structures
$C^\infty([0;1]^n\times\mathbb{R}^n)$ and $compl_{{\mathcal
T}{\mathcal G}_0} {\mathcal T}{\mathcal C}$ respectively. The
chain $\Phi \in Ch_n(M)$ is \emph{smoothly extendable with respect
to} ${\mathcal G}$ when every cube $\varphi \in supp \Phi$ is
smoothly extendable with respect to ${\mathcal G}$. The smooth
point $k$-form $\eta$ on $M$ is said to be \emph{smoothly
extendable with respect to} ${\mathcal G}$ if there exists a
smooth point $k$-form $\tilde{\eta}$ on $compl_{\mathcal G} M$,
such that $\eta = \iota^*(\tilde{\eta})$, where $\iota : M
\rightarrow compl_{\mathcal G} M$ is the natural embedding (for a
point $p \in M$ we put its filter of neighborhoods $\iota(p) =
{\mathcal F}_p$ --
see \cite{cuk}).\\
\end{df}

\smallskip
Let $\pi : TM \rightarrow M$ be the canonical projection. Then for
each set of generators ${\mathcal G}$ of the differential
structure ${\mathcal C}$ the mapping $\pi$ is uniform with respect
to uniform structures ${\mathcal U}_{{\mathcal T}{\mathcal G}_0}$
and ${\mathcal U}_{\mathcal G}$ (from definitions of $\pi$ and
${\mathcal U}_{{\mathcal T}{\mathcal G}_0}$). Moreover for each
differential space $(N,{\mathcal D})$ and each $F : (N,{\mathcal
D}) \rightarrow (M, {\mathcal C})$ we have $\pi \circ TF = F$. So
the fact that a smooth generalized $n$-dimensional cube $\varphi :
{\mathcal D} \rightarrow M$ is smoothly extendable  implies that
$\varphi$ is a uniform mapping and its extension $\tilde{\varphi}
: [0,1]^n \rightarrow compl_{\mathcal G} M$ is smooth, i.e. it is
a smooth $n$-dimensional cube on $compl_{\mathcal G} M$.

\begin{df}
Let $\varphi : D \rightarrow M$ be a smooth generalized
$n$-dimensional cube on a differential space $(M,{\mathcal C})$,
smoothly extendable with respect to some family of generators
${\mathcal G}$ of the structure ${\mathcal C}$. We define
\emph{the boundary of $\varphi$ with respect to ${\mathcal G}$} as
the $n$-dimensional smooth chain $\partial_{\mathcal G} \varphi$
on $compl_{\mathcal G} M$ given by the formula:
$$\partial_{\mathcal G} \varphi = \sum_{i=1}^n \sum_{\alpha=
0,1}(-1)^{i + \alpha}\tilde{\varphi}_{(i,\alpha)},$$ where
$\tilde{\varphi}_{(i,\alpha)} = \tilde{\varphi} \circ
(I_{(i,\alpha)}^n)$ is \emph{the $(i,\alpha)$-face} of the cube
$\tilde{\varphi}$,
$$I_{(i,0)}^n (x) = (x^1,\ldots,x^{i-1},0,x^i,\ldots,x^{n-1}),$$
$$I_{(i,1)}^n (x) = (x^1,\ldots,x^{i-1},1,x^i,\ldots,x^{n-1})$$
and $\tilde{\varphi}$ is the smooth extension of $\varphi$ on
$compl_{\mathcal G} M$. \label{67}
\end{df}

It is easy to see that $\partial_{\mathcal G}\varphi =
\partial_{\tilde{\mathcal G}} \tilde{\varphi}$, where
$\tilde{{\mathcal G}}$ is the family of extensions of elements of
the set ${\mathcal G}$ on $compl_{\mathcal G} M$.

\begin{df}
\emph{The boundary of a generalized n-dimensional chain} $\Phi =
\sum_{j=1}^m c_j \varphi_j$ on the differential space
$(M,{\mathcal C})$ smoothly extendable with respect to a family of
generators ${\mathcal G}$ of the structure ${\mathcal C}$ we call
the $(n-1)$-dimensional chain:
$$\partial_{\mathcal G} \Phi := \sum_{j=1}^m c_j \partial_{\mathcal G} \varphi_j,$$
where $\{\varphi_1,\ldots,\varphi_m\} = supp \Phi$.
\end{df}

For formulating the analogue of Stokes theorem we need a mapping
that will be an analogue of the exterior derivative. That mapping
should be a generalization of the exterior derivative on manifolds
and it should linearly maps $k$-forms to $(k+1)$-forms. It should
also commute with pullback of forms. Next example shows that we
can not give such mapping uniquely even on smooth extendable
forms.

\begin{prz}
Let $M = \{(x,y) \in \mathbb{R}^2 : xy = 0\}$, and ${\mathcal C} =
C^\infty (\mathbb{R}^2)_M$. Then the zero $1$-form on
$(M,{\mathcal C})$ could be extended to the zero $1$-form on
$\mathbb{R}^2$ and to the form $\eta = x dy - y dy$. We have:
$d\eta = 2 dx \wedge dy$. Moreover if by $\iota : M \rightarrow
\mathbb{R}^2$ we denote the natural embedding, then $\iota^* (dx
\wedge dy) = dx \wedge dy \neq 0$, because $dx \wedge dy_{|T_0 M}
\neq 0$, while $\iota^*(0) = 0$. $\qquad \Box$
\end{prz}

If ${\mathcal G}$ is a family of bounded generators of a
differential structure ${\mathcal C}$ on a set $M$ (we can assume
that if $g \in {\mathcal G}$, then $|g| \leq 1$), then the
completion $compl_{\mathcal G} M$ is a differential
compactification $compt_{\mathcal G} M$ (see \cite{cuk}).

\begin{stw}
Let the previous assumption for ${\mathcal G}$ be satisfied. Then
any point $k$-form $\eta$ that is smoothly extendable with respect
to ${\mathcal G}$ can be extended to a smooth $k$-form $\eta _1$
on $[-1,1]^{\mathcal G}$.
\end{stw}

\emph{Proof}. Let $\tilde{\eta}$ be such a smooth point $k$-form
on $compl_{\mathcal G} M$ that $\eta=i^*\tilde{\eta}$, where $i$
is the natural embedding $M$ into $compl_{\mathcal G} M$. Let
$\tilde{{\mathcal G}}$ be the set of all continuous extensions of
elements of ${\mathcal G}$ onto $compl_{\mathcal G} M$. Since
$\tilde{{\mathcal G}}$ is a family of generators of the
differential structure $compl_{\mathcal G} {\mathcal
C}=C^\infty([-1;1]^{\mathcal G})_{compl_{\mathcal G} M}$ we can
choose for any $p\in compl_{\mathcal G} M$ a neighborhood $W_p$
such that (\ref{4}) holds for the form $\tilde{\eta}$ and for some
$\tilde{\alpha}_{i_j}\in\tilde{\mathcal G}$. If we take a set
$V_p$ open in $[-1;1]^{\mathcal G}$ and such that
$W_p=compl_{\mathcal G}M \cap V_p$ then the same formula defines
the point $k$- form $\tilde{\eta}_p$ on $V_p$. But the space
$compl_{\mathcal G} M$ is compact and therefore there exists
$n\in\mathbb{N}$ and $p_1,\ldots,p_n\in compl_{\mathcal G} M$ such
that $compl_{\mathcal G} M\subset\bigcup\limits_{j=1}^nV_{p_j}$.
Putting $V_0:=[-1;1]^{\mathcal G}\setminus compl_{\mathcal G} M$,
$V_j:=V_{p_j}$ and $\tilde{\eta}_j:=\tilde{\eta}_{p_j}$ for
$j=1,\ldots,n$ we obtain the open finite covering
$\{V_0,V_1,\ldots,V_n\}$ of $[-1;1]^{\mathcal G}$. Let
$\{\gamma_0, \gamma_1,\ldots,\gamma_n\}$ be a smooth partition of
unity subordinated to the covering. Than any
$\gamma_j\tilde{\eta}_j,\ \ j=1,\ldots,n$, can be in a standard
manner extend to $[-1;1]^{\mathcal G}$ (we take
$(\gamma_j\tilde{\eta}_j)_{|[-1;1]^{\mathcal G}\setminus W_j}=0)$.
Taking $\eta_1:=\sum\limits_{j=1}^n \gamma_j\tilde{\eta}_j$ we
obtain the $k$-form on $[-1;1]^{\mathcal G}$ which fulfils the
conditions of the proposition. \hfill$\Box$

\bigskip
On the differential space $([-1,1]^{\mathcal G}, C^\infty
([-1,1]^{\mathcal G}))$ there is the well definite operator of
exterior derivative $d$ which assigns to any smooth point $k$-form
$\omega$ the smooth $(k+1)$-form $d\omega$ given by the formula
$$(d \omega)(v_1,\ldots,v_{k+1}) = \sum_{i=1}^{k+1} (-1)^{i+1}
X_i (\omega (X_1,\ldots,\hat{X}_i,\ldots,X_{k+1}))(p) +$$
\begin{equation}+ \ldots \sum_{i<j}(-1)^{i+j} \omega ([X_i,X_j],X_1,\ldots,
\hat{X}_i,\ldots,\hat{X}_j,\ldots,X_{k+1})(p),
\label{89}\end{equation} where $v_1,\ldots,v_{k+1} \in
T_p([-1,1]^\mathcal{G})$, $X_1,\ldots,X_{k+1}$ are vector fields
on $[-1,1]^\mathcal{G}$, such that $X_j(p) = v_j$ for
$j=1,2,\ldots,n$ (such a vector fields exist on a space which is a
Cartesian product of intervals; we can take $X_j=const=v_j\in
T[-1;1]^{\mathcal G}\cong\mathbb{R}^{\mathcal G}$ -- see
\cite{DDD}, \cite{cuk}, \cite{DPW1} or \cite{DPW2}) and
$"^\wedge"$ means that the variable disappears in that expression.

For any smooth mapping $F:[0;1]^n \rightarrow [-1,1]^{\mathcal G}$
and any smooth point $k$-form $\omega$ smoothly extendable with
respect to ${\mathcal G}$ operator $d$ fulfils the equality:
$d_1(F^* \chi) = F^* d\chi$, where on the left hand side there is
"the classical" operator $d_1$ of the exterior derivative on
$[0;1]^n$. It is the result of (\ref{89}), definitions and
properties of vector fields. Now we can get the analogue of Stokes
theorem.

\begin{tw}
If a smooth point $k$-form $\eta$ and a smooth generalized
n-dimensional chain $\Phi$ are smoothly extendable with respect to
a family ${\mathcal G}$ of bounded generators of a differential
structure ${\mathcal C}$ on a set $M$, then
\begin{equation}
\int_{\partial_{\mathcal G} \Phi} \eta = \int_{\Phi} d
\tilde{\eta}, \label{90}\end{equation} where $\tilde{\eta}$ is an
arbitrary smooth extension of the $k$-form $\eta$ on the
differential space $([-1;1]^{\mathcal G},C^\infty([-1;1]^{\mathcal
G}))$ and $d$ is the operator of exterior derivative given by
(\ref{89}). \label{43}
\end{tw}

\emph{Proof}. It is enough to show equality (\ref{90}) for any
$n$-dimensional cube $\varphi$ on $(M,{\mathcal C})$ smoothly
extendable with respect to ${\mathcal G}$. We have
$$\int_{\partial_{\mathcal G} \varphi} \eta = \int_{\partial([0;1]^n)}
\widetilde{\varphi^* \eta} = \int_{\partial ([0;1]^n)}
\tilde{\varphi}^* \tilde{\eta} = \int_{[0;1]^n} d_1
(\tilde{\varphi}^* \tilde{\eta})= \int_{[0;1]^n} \tilde{\varphi}^*
d\tilde{\eta} = \int_\varphi d\tilde{\eta},$$ where $d_1$ is the
operator of exterior derivative on $[0;1]^n$ and $\tilde{\varphi}$
is the extension of $\varphi$ on $[0;1]^n$. \hfill $\Box$\\

\vspace{1cm}

\footnotesize
\begin{center}

\end{center}

\end{document}